\documentclass[12pt]{article}
\usepackage{amsmath}
\usepackage{amsfonts}
\usepackage{amscd}
\newcommand{\qed}{\hfill \fbox{}\medskip}

\begin{document}
\title{Subalgebras of the Gerstenhaber algebra of 
differential operators}
\author{A.V.Bratchikov\\
 Department of Mathematics\\ Kuban
State Technological University,\\ 
Krasnodar, 350072, Russia 
}
 \date{} \maketitle
\begin{abstract}

We construct a family of subalgebras of the Gerstenhaber algebra of differential operators. The subalgebras are labeled by subsets of the additive group ${\mathbb Z}^n$ that are closed under addition.   Each subalgebra is invariant under the Hochschild coboundary operator. 
\end{abstract}


The Gerstenhaber algebra of cochains of a commutative algebra was introduced in \cite{G}. It is equipped with two binary operations: a cup product $\smile$ and a bracket $[\,\cdot\,,\cdot\,].$ Hochschild cocycles and coboundaries form subalgebras of the Gerstenhaber algebra. Lie algebras of vector fields are subalgebras of the Gerstenhaber bracket algebra of differential operators. These algebras play an important role in deformation quantization of Poisson manifolds \cite{BFFLS},\cite{PV},\cite{K}.

In this paper a family of subalgebras of the Gerstenhaber algebra of differential operators
is constructed. We find a commutative subalgebra 
that splits 
the Gerstenhaber algebra 
into a direct sum of weight subspaces.
The weight
subspaces are labeled 
by elements of the additive group $\mathbb Z^n.$
The subalgebras are constructed as direct sums of weight subspaces. To each additive semigroup $\Delta\subset \mathbb Z^n$ 
we put into correspondence a differential Gerstenhaber algebra $C_\Delta$ with the Hochshild differential.
Elements of $\Delta$ define a grading of $C_\Delta.$ 

As an application of these results we consider deformation quantization of ${\mathbb R^2}$. By using the Maurer-Cartan equation we show that if a Poisson bivector belongs to an algebra $(C_\Delta, [\,\cdot\,,\cdot\,]),$ then higher order cochains of the corresponding star product can be chosen from an ideal of this algebra.

Let 
$A$ be a space of real 
polynomials in $x_i,$ $i=1,\ldots,n,$ ${A=\mathbb R[x_1
,\ldots,x_n
],}$ and let $ 
C^p= Hom(A^p, A),$ $p\ge 0.$ $C^p$
consists of the operators of the form 
\begin{eqnarray*}
\varphi =\sum_{a_0a_1\ldots a_p}
\varphi_{a_0a_1\ldots a_p}x^{a_0}\partial^{a_1}
\otimes \ldots \otimes \partial^{a_p},
\end{eqnarray*}
where $\varphi_{a_0a_1\ldots a_p}
$
are real constants, $x^{a}=x^{a^1}_1x^{a^2}_2\ldots x^{a^n}_n,$ 
${a=(a^1,a^2,\ldots,{a^n})
\in 
\mathbb N
^n,}$ $\mathbb N$ is the set of nonnegative integers, $
\partial^{a}=
\partial_1^{a^1}\partial_2^{a^2}\ldots \partial^{a^n}_n,$ ${
\partial_{i}=
\partial/ \partial x_i.}$ 
For $u_1,u_2,\ldots, u_p\in A$
\begin{eqnarray*}\varphi( u_1,u_2,\ldots, u_p)= \sum_{ a_0a_1\ldots a_p}
\varphi_{a_0a_1\ldots a_p}z^{a_0}
(\partial^{a_1}u_1)\ldots  (\partial^{a_p}u_p).
\end{eqnarray*}
We set $
 C= \bigoplus_{p\ge 0}
 C^p.$

The cup product and the Gerstenhaber bracket 
on $C$ are given by
\begin{eqnarray}\label{C}\varphi\smile \psi (u_1,\ldots, u_{p+q}) = \varphi(u_1,\ldots  u_p)\psi(u_{p+1}, \ldots u_{p+q}),
\end{eqnarray}
\begin {multline} \label{U}
[\varphi, \psi](u_1,\ldots,u_{p+q-1})=\\ 
=\sum_{k=1}^p(-1)^{(k-1)(q-1)}
\varphi (u_1,\ldots,u_{k-1},\psi(u_{k},\ldots, u_{k+q-1}), u_{k+n} ,\ldots ,u_{p+q-1})-\\
-(-1)^{(p-1)(q-1)}\sum_{k=1}^q(-1)^{(k-1)(p-1)}
\psi (u_1,\ldots,u_{k-1},\varphi(u_k,\ldots, u_{k+p-1}), u_{k+p} ,\ldots ,u_{p+q-1}),
\end {multline}
for $\varphi\in C^p,$ $\psi\in C^q.$ 
The cup product is associative. The bracket is a skew multiplication satisfying the Jacobi identity, i.e. if $\varphi\in C^p,$ $\psi\in C^q,$  $\chi\in C^r,$ then  
$$[\varphi,\psi]=-(-1)^{(p+1)(q+1)}[\psi,\varphi],$$
\begin{eqnarray} 
\label {J} (-1)^{(p+1)(r+1)}[\varphi,[\psi,\chi]]+\mbox{cycl. perm.}\, (\varphi,\psi,\chi)
=0.
\end{eqnarray}

Let ${\bf \delta}$ denote the Hochschild coboundary operator. 
It is defined by
\begin {multline} 
\label {nut} 
\delta \varphi (u_1,u_2,\ldots, u_{p+1})= u_1\varphi (u_2,\ldots, u_{u+1})+ \\
+\sum_{k=1}^p(-1)^k
\varphi (u_1,\ldots, u_{k-1},u_k u_{k+1},u_{k+2},\ldots, u_{p+1})+(-1)^{p+1}\varphi (u_1,\ldots, u_p)u_{p+1}.
\end {multline}
The differential $\delta$ can be written in terms of the bracket  $[\,\cdot\,,\cdot\,]$ as 
\begin{eqnarray} 
\label {Jo}\delta \varphi =-[\varphi,m],
\end{eqnarray}
where $m(u_1,u_2)$ is the product $ u_1u_2$ of elements $u_1,u_2\in A.$

An algebra $({\cal A}, [\,\cdot\,,\cdot\,],\delta)$ is called a differential 
Gerstenhaber algebra \cite {HS}. 
In this paper by a differential 
Gerstenhaber algebra we shall mean an algebra ${({\cal A},\smile,[\,\cdot\,,\cdot\,],\delta).}$ 

{\bf Lemma 1.} {\it If $\chi$ is a vector field, $\chi=\sum_
i
\chi^i(x)\partial_i,$ $\varphi\in  C^p,$ $\psi\in C^q,$    then
\begin{eqnarray}\label{B}[\chi,\varphi\smile\psi]= [\chi,\varphi]\smile\psi+\varphi\smile [\chi,\psi].
\end{eqnarray}
}
{\bf Proof.} We have 
\begin {multline*} 
[\chi,\varphi\smile\psi](u_1,\ldots,u_{p+q})=\\ 
=
\bigl(\chi\varphi (u_1,\ldots,u_{p})\bigr)\psi(u_{p+1},\ldots, u_{p+q})+\varphi (u_1,\ldots,u_{p})\bigl(\chi\psi(u_{p+1},\ldots, u_{p+q})\bigr)-\\
-\sum_{k=1}^{q+p}
\varphi (u_1,\ldots,u_{k-1},\chi(u_k),u_{k+1},\ldots, u_{p})\psi (u_{p+1} ,\ldots ,u_{p+q})=\\
=\bigl([\chi,\varphi]\smile\psi\bigr)(u_1,\ldots,u_{p+q})+\bigl(\varphi\smile [\chi,\psi]\bigr)(u_1,\ldots,u_{p+q}). 
\end {multline*}

Vector fields $h^i =x_i\partial_i,$ $i=1,\ldots n,$ generate a commutative subalgebra of $(C, [\,\cdot\,,\cdot\,]), $ 
$$
[h^i,h^j]=0.
$$
From (\ref{U}) we get
$$
[h^i, 
x^{a_0}\partial^{a_1}\otimes\partial^{a_2}\otimes \ldots \otimes \partial^{a_p}]=
({a_0^i}-\sum_{s=1}^p{a_s^i})x^{a_0}\partial^{a_1}\otimes \partial^{a_2}\otimes \ldots \otimes \partial^{a_p}.
$$
Let $ C_a$ be a subspace of $C,$ spanned by the cochains   
\begin{eqnarray*}
x^{a_0}\partial^{a_1}\otimes \partial^{a_2}\otimes \ldots \otimes \partial^{a_p}
\end{eqnarray*}
with $
{a_0}-\sum_{s=1}^p
{a_s}=a,$ $p\ge 0.$

{\bf Theorem 1.}
 {\it Let  
 $C_\Delta = \bigoplus_{a\in \Delta} C_a,$ where $\Delta \subset \mathbb Z^n$ is an additive semigroup.
Then $
(C_\Delta, \smile, [\,\cdot\,,\cdot\,],\delta)
$ is a differential Gerstenhaber algebra.}

This result follows from the following Lemma.

{\bf Lemma 2.} {\it If $\varphi\in  C_a,$ $\psi\in C_b,$ then
$$\varphi\smile\psi\in  C_{a+b},\qquad [\varphi,\psi]\in C_{a+b},\qquad \delta\varphi\in  C_{a}.$$}
{\bf Proof.}The first two statements 
follow from  identities (\ref{J}) and (\ref{B}) with $\chi=h^i.$ For $\varphi\in  C_{a}$ 
(\ref {Jo}) implies $\delta\varphi\in  C_{a},$ since $m\in C_0.$
$ \qed$

{\bf Proposition 1.} {\it (i) For any $r\ge 1$ $ I^{(r)}_\Delta=\bigoplus_{a_1\ldots a_r\in \Delta}C_{a_1+\ldots +a_r}$ is an ideal of $C_\Delta.$

(ii) $ I^{(r)}_\Delta=C_\Delta$ if and only if $ \Delta+\Delta=\Delta$.}

{\bf Proof.} (i) follows from Lemma 2. 
In order to prove (ii),   
we observe that $ \underbrace{\Delta+\ldots +\Delta}_r = \Delta$ if and only if $\Delta+\Delta=\Delta.\qed$

Proposition 1 states, in particular, that if $\Delta+\Delta\ne \Delta,$ then $I^{(r)}_\Delta,$ $r\ge 2,$ is a nontrivial ideal of $C_\Delta.$
 
Define maps $\theta_i:C\to C
,$ $i=1,\ldots, n,$ by
$$
\theta_i \varphi
=(-1)^{
a^i}\varphi
$$
for $\varphi\in C_a.$
For a multi-index $I = (i_1, \ldots, i_r),$ $
\theta_{I}$ will stand for $\theta_{i_1}\theta_{i_2}\ldots \theta_{i_r}.$ 
From the relation $C= \bigoplus_{a\in \mathbb Z^n} C_a$ and Lemma 2 it follows that for any ${1\le r\le n,}$
$
\theta_{I}$ 
is an involutive automorphism of $(C, \smile, [\,\cdot\,,\cdot\,]).$ 
Let 
\begin{eqnarray}\label{inv}
 C_{I}^+=\{\varphi\in C\,|\,\theta_{I}\varphi=\varphi\},\qquad C_{I}^-=\{\varphi\in C\,|\,\theta_{I}\varphi=-\varphi\}.
\end{eqnarray}
The space $C_{I}^+$ can be represented as 
$$
 C_{I}^+=\bigoplus_{a\in \Theta_{I}}C_a,
$$
where $\Theta_{I}$ is a subgroup of $\mathbb Z^n,$ defined by
$$\Theta_{I}=\{a\in\mathbb Z^n\,|\,a^{i_1}+a^{i_2}+\ldots+a^{ i_r}\,\, \mbox{even}\}.$$
Then $\delta C_{I}^+\subset C_{I}^+,$ by Lemma 2.

Let $g$ be an algebra with an involutive automorphism $\theta.$ Recall that the subset  
$$l=\{\upsilon \in g\,|\,\theta \upsilon=\upsilon\}$$   
is a subalgebra which is called an involutive algebra (of involutive automorphism $\theta$)
\cite{S}. 
The following is an immediate consequence of these results. 

{\bf Proposition 2.} {\it $ C_{I}^+
$ is an involutive differential Gerstenhaber algebra.
}

By using (\ref{inv}), we get
$$ C= C_{I}^+\oplus  C_{I}^-,$$ 
$$ C_{I}^+\smile\,C_{I}^+\subset C_{I}^+,\qquad C_{I}^+\smile C_{I}^-\subset  C_{I}^-,
$$
$$
C_{I}^-\smile C_{I}^+\subset  C_{I}^- ,\qquad C_{I}^{- }\smile C_{I}^-\subset C_{I}^+,
$$
$$
[ C_{I}^+,  C_{I}^+]\subset C_{I}^+,\qquad
[ C_{I}^+, C_{I}^-]\subset C_{I}^-,\qquad
[ C_{I}^-, C_{I}^-]\subset  C_{I}^+.$$

{\bf Proposition 3.} {\it Let $ C_H= \bigoplus_{a\in H} C_a,$ $H\subset \mathbb Z^n,$ and let ${C =C_H \oplus L.}$ Then $$C_H\smile L\subset L,\qquad
L\smile C_H\subset L,\qquad
[C_H, L]\subset  L$$
if and only if $H$ is a subgroup of $\mathbb Z^n.$}

{\bf Proof.}
The space $L$ is given by 
$L =\bigoplus_{a\in K}C_a,$ $K={\mathbb Z}^n\setminus H.$
It is easily verified that 
the relation 
$H+K \subset K$ holds if and only if $H$ is a subgroup of $\mathbb Z^n.$  $\qed$

Denote by $ C_{ab}$ a subspace of $ C,$ spanned  by the operators   
$$x^{a_0}\partial^{a_1}\otimes \partial^{a_2}\otimes \ldots \otimes \partial^{a_p},
\qquad 
p\ge 0,
$$
with $
{a_0}-\sum_{s=1}^p
{a_s}=a
,$ $
{a_0}+\sum_{s=1}^p
{a_s}=b
.$ Let $\le$ be a lexicographical ordering on $\mathbb Z^n.$ Then $C_a=\bigoplus_{a\le b}C_{ab}.$
A lexicographical ordering on $\mathbb Z^n\times \mathbb Z^n$ is defined as $(a,b)\le (a',b')$ if and only if $a<a'$ or $(a=a'\mbox{ and } b\le b').$ Let  $(J,\le)$ be an additive ordered semigroup that consists of $(a,b)\in \mathbb Z^n\times \mathbb Z^n,$ 
$a\le b.$ Define a space $S_\alpha$ by
$$S_\alpha=S_{ab}=\bigoplus_{a\le d\le b}C_{ad}
,\qquad \alpha =(a,b)
.$$

{\bf Proposition 4.}
{\it The set $(S_\alpha,$ $\alpha\in J)$ is a filtration of $
(C,\smile, [\,\cdot\,,\cdot\,],\delta)
.$}

{\bf Proof.}
From (\ref{C}), (\ref{U}) and Lemma 2 it follows that
$$S_{\alpha}\smile S_{\beta}\subset S_{\alpha+\beta},\qquad [S_{\alpha}, S_{\beta}]\subset S_{\alpha+\beta}, \qquad \delta S_{\alpha}\subset S_{\alpha}.$$ It is easily seen that
$$ {C= \bigcup_{\alpha\in J}S_{\alpha},\qquad S_{\alpha}\subset S_{\beta},\qquad {\alpha}< {\beta}
. }\qquad \phantom{x}
$$
Hence the result. $\qed$

This filtration can be used to define an algebra $(S,\smile, [\,\cdot\,,\cdot\,],\delta),$ that 
consists of the elements of the form
$$\sum_{\alpha \in J} \tau^\alpha v_{\alpha}, 
$$
where $\tau^{\alpha}=r^{a}s^b,$  $r^{a}=r^{a^1}_1r^{a^2}_2\ldots r^{a^n}_n,$ 
$s^{b}=s^{b^1}_1s^{b^2}_2\ldots s^{b^n}_n,$  $(r_i,s_i,$ $i=1,\ldots,n)$ are indeterminates,  $v_{\alpha}\in S_{\alpha}.$

A star product on $\mathbb R^n$ is given by
\begin{eqnarray*}
f\ast g=fg+\pi(f,g),
\end{eqnarray*}
where  
the 2-cochain $\pi=(t/ 2) \pi_1+ O(t^2),$  $t\in {\mathbb R},$ satisfies the 
Maurer-Cartan equation 
\begin{eqnarray} \label {maurer}
\delta \pi=\frac 1 2 [\pi,\pi].
\end{eqnarray} 
Here $\pi_1$ is necessarily a Poisson bivector. 
Let $C_\Delta$ be as in Theorem 1.

{\bf Proposition 5.} {\it
Let
$\pi=
(t/ 2) \pi_1+ \sum_{k\ge 2}^\infty t^k
\pi_k$ be a solution to (\ref{maurer}) in the case of ${\mathbb R^2}$.  
If ${\pi_1 \in C_\Delta,}$ 
then $\pi_k,$ $k\ge 2,$ can be chosen from $I^{(2)}_\Delta.$}

{\bf Proof.} Eq. (\ref{maurer}) can be expressed by an infinite set of relations
\begin{eqnarray} 
\label {maurer1}
\delta \pi_k=B_k,
\qquad k=1,2,\ldots,
\end{eqnarray}
where 
\begin{eqnarray} 
\label {maurer11}
B_k=\frac 1 2\sum_{i+j=k} [\pi_i,\pi_j].\qquad k=1,2,\ldots
\end{eqnarray} 
Assume that $\pi_i\in I^{(2)}_\Delta,$ $i\le k-1,$ have been constructed. It is directly verified that $\delta B_k=0$ \cite{PV}. Therefore, in the case of ${\mathbb R^2}$ the $3$-cochain $B_k$ is a coboundary, 
$B_k=\delta X_k$ \cite{GR}. This means that (\ref{maurer1}) possesses a solution.
Eq. (\ref{maurer11}) implies $B_k\in I^{(2)}_\Delta.$ From this and Lemma 2 it follows that 
$\pi_k=X_k+ \Upsilon_k,$ where  $X_k\in I^{(2)}_\Delta,$
$\Upsilon_k \in Ker\, \delta.$ On ${\mathbb R^2}$ one can put $\Upsilon_k=0.$ $\qed$

\end{document}